\documentclass[12pt]{revtex4}

\usepackage{amssymb,color}
\usepackage{algorithmicx, algpseudocode}
\usepackage{hyperref}

\begin{document}

\title{Geodesic acceleration and the small-curvature approximation for nonlinear least squares}

\author{Mark K.~Transtrum}
\affiliation{Department of Bioinformatics and Computational Biology, University of Texas M.D.~Anderson Cancer Center, Houston Texas, U.S.A.}
\email{mkt26@cornell.edu}

\author{James P.~Sethna}
\affiliation{Laboratory of Atomic and Solid state Physics, Cornell University, Ithaca, New York 14853, USA}

\begin{abstract}

It has been shown numerically that the performance of the Levenberg-Marquardt algorithm can be improved by including a second order correction known as the geodesic acceleration.  In this paper we give the method a more sound theoretical foundation by deriving the geodesic acceleration correction without using differential geometry and showing that the traditional convergence proofs can be adapted to incorporate geodesic acceleration.  Unlike other methods which include second derivative information, the geodesic acceleration does not attempt to improve the Gauss-Newton approximate Hessian, but rather is an extension of the small-residual approximation to cubic order.  In deriving geodesic acceleration, we note that the small-residual approximation is complemented by a small-curvature approximation.  This latter approximation provides a much broader justification for the Gauss-Newton approximate Hessian and Levenberg-Marquardt algorithm.  In particular, it is justifiable even if the best fit residuals are large, is dependent only on the model and not on the data being fit, and is applicable for the entire course of the algorithm and not just the region near the minimum.

\end{abstract}
\maketitle
\section{Introduction}
\label{sec:introduction}

In this paper we consider the problem of minimizing a scalar function whose form is the sum of squares
\begin{equation}
  \label{eq:Cost}
  C(\theta) = \frac{1}{2} \sum_m r_m(\theta)^2,
\end{equation}
where $r: \mathbf{R}^N \rightarrow \mathbf{R}^M$ is an $M$-dimensional vector function of $N$ parameters, $\theta$.  We refer to $r(\theta)$ as the residuals and scalar function $C(\theta)$ as the cost.  Functions of this form often arise in the context of data fitting and represent an important class of problems as evidenced by the large number of software packages dedicated to their optimization.

The structure of this particular problem lends itself to efficient optimization.  In particular, consider the Hessian matrix of second derivatives, necessary to implement a quasi Newton method:
\begin{eqnarray}
  \label{eq:Hessian}
\frac{\partial^2 C }{\partial \theta_\mu \partial \theta_\nu} & = & \sum_m \left( \frac{\partial r_m}{\partial \theta_\mu} \frac{\partial r_m}{\partial \theta_\nu} + r_m \frac{\partial^2 r_m}{\partial \theta_\mu \partial \theta_\nu} \right)  \\
 & \approx & \sum_m \frac{\partial r_m}{\partial \theta_\mu} \frac{\partial r_m}{\partial \theta_\nu} \nonumber \\
 & = & (J^TJ)_{\mu\nu}, \nonumber
\end{eqnarray}
where in the last two lines we have applied the so-called Gauss-Newton  or small-residual approximation and introduced the Jacobian matrix of first derivatives $J$.  The approximation in the second line is usually justified by the hope that at a minimum of $C(\theta)$, the individual residuals are small, so that the Hessian is dominated by the contributions from the first term.  Numerically, this approximation is advantageous since it allows one to implement a quasi-Newton method by calculating only the first derivatives of the residuals.  This approximation comes at the cost of storing the derivative information for each residual individually, but this is rarely a bottleneck on modern computers.  Consequently, the functional form of Eq.~(\ref{eq:Cost}) effectively allows one to estimate the Hessian matrix with the same information used to calculate the gradient.  

Applying a trust region method to the Gauss-Newton approximate Hessian results in the Levenberg-Marquardt algorithm\cite{Levenberg1944,Marquardt1963,Osborne1976,More1977,Nocedal2000,Press2007} which iteratively updates the parameters according to
\begin{equation}
  \label{eq:LMstep}
  \delta \theta = - \left( J^TJ + \lambda D^TD \right)^{-1} g,
\end{equation}
where $\lambda$ is an appropriately chosen Langrange multiplier for the step bound $\delta \theta^T D^T D \delta \theta \leq \Delta^2$ and $g = J^T r$ is the gradient.  Because Levenberg-Marquardt is a quasi-Newton method, it usually has very good convergence properties.  In particular, if the small-residual approximation is good, convergence can be super-linear to a local minimum.  Furthermore, for well-constructed choices of $\lambda$ and $\Delta$, the method is globally convergent\cite{Osborne1976,Nocedal2000}.  

Although the Levenberg-Marquardt algorithm has many desirable properties, it is not always ideal.  Often data fitting problems have a cost function characterized by narrow, winding canyons.  Although, asymptotically the algorithm may converge super-linearly, it may nevertheless spend an unreasonable amount of time navigating the winding canyon before it finally zooms into the minimum.  This problem is typically more severe on problems with many parameters, which in turn are often more computationally expensive to evaluate and lack good parameter estimates to use as starting points.  An improved optimization method which can find minima with fewer function evaluations (and especially fewer Jacobian evaluations) would be a welcome improvement.

In order to help improve the Levenberg-Marquardt algorithm, the authors previously proposed the inclusion of a geodesic acceleration term in the algorithm\cite{Transtrum2010,Transtrum2011}.  This correction was derived using an information geometric interpretation of the least-square problem and was justified based on the empirically observed small-extrinsic curvatures of the relevant manifolds.  In this paper, we will see that the geodesic acceleration correction can be understood as a generalization of the Gauss-Newton method extended to cubic order.  Although other methods exist which utilize higher-order information, geodesic acceleration is complementary to these approaches as their primary motivation is to improve the estimate of the Hessian.  By contrast, the geodesic acceleration assumes the Hessian estimate is adequate to proceed to higher-order.

In this paper, we derive the geodesic acceleration in a geometric independent way (section \ref{sec:derivation}) and prove that its inclusion in the Levenberg-Marquardt algorithm does not compromise its convergence properties (section \ref{sec:Convergence}).  From the explicit derivation in section \ref{sec:derivation}, we see that geodesic acceleration can be understood as a continuation of the small-residual approximation to higher-order terms.  Indeed, the discarded terms are properly understood as the residuals coupled to the extrinsic curvatures of the \emph{Model Graph} defined in references\cite{Transtrum2010,Transtrum2011}.  The small-curvature approximation, therefore, provides additional justification for the Gauss-Newton approximation and the geodesic acceleration correction.  As we argue in section \ref{sec:smallcurvature}, the small-curvature approximation is more useful than the small-residual approximation for many reasons.  In particular, it is valid not only near a local minimum of the cost, but for all parameter values; it is a property of the model and not the data being fit and so is valid even when the model cannot fit the data well; furthermore, numerical experiments on many models suggest the small-curvature approximation is nearly universally valid.

\section{Derivation}
\label{sec:derivation}

In order to improve the efficiency of the Levenberg-Marquardt method, we propose modifying the step to include higher order corrections in a numerically efficient manner.  To derive this correction, consider the minimization problem of finding the best residuals with a constrained step-size.  We write the dependence of the residual on the shift $\delta \theta$ as
\begin{equation}
r(\theta + \delta\theta) =  r + J \delta\theta 
	+ 1/2\, \delta\theta^T K \delta\theta + \cdots,
\end{equation}
where $J$ and $K$ are the arrays of first and second derivatives respectively.  We wish to minimize
\begin{equation}
  \label{eq:minnolam}
  \min_{\delta\theta} \ \ \left( r + J \delta\theta 
	+ 1/2\, \delta \theta^T K \delta \theta \right)^2
\end{equation}
with the constraint that $\delta \theta^T D^T D \delta \theta \leq \Delta^2$.  After introducing a Lagrange multiplier $\lambda$ for the constraint in the step size, the minimization becomes
\begin{equation}
  \label{eq:minlam}
  \min_{\delta \theta}  \left( r + J \delta \theta
	+ 1/2\, \delta \theta^T K \delta \theta \right)^2 
	+ \lambda \delta \theta^T D^T D \delta \theta.
\end{equation}

By varying $\delta \theta$ we find the normal equations:
\begin{eqnarray}
  \label{LM:eq:normal}
  \sum_m J_{m\mu} r_m + \sum_{m\nu} \left( J_{m\mu} J_{m\nu} + r_m K_{m\mu\nu} + \lambda D_{m\mu} D_{m\nu} \right) \delta \theta_\nu & &\nonumber \\
  + \sum_{m\nu\alpha} \left( J_{m\nu} K_{m\mu\alpha} + 1/2\, J_{m\mu} K_{m\nu\alpha} \right) \delta \theta_\nu \delta \theta_\alpha & =& 0,
\end{eqnarray}
where we have explicitly included all the indices to avoid any ambiguity.  Since we constrain the step size, it is natural to assume that $\delta \theta$ is small, and we seek a solution of Eq.~(\ref{LM:eq:normal}) as a perturbation series around the linearized equation:
\begin{equation}
  \label{eq:stepseries}
  \delta \theta = \delta \theta_1 + \delta \theta_2 + \cdots.
\end{equation}
Let $\delta \theta_1$ be a solution of the linearized equation:
\begin{eqnarray}
  \delta \theta_1 & = & -(J^TJ + r^T K + \lambda D^TD)^{-1} J^T r \nonumber\\
  & \approx & -(J^TJ + \lambda D^TD)^{-1} J^T r, \nonumber
\end{eqnarray}
where in the second line we have made the usual Gauss-Newton approximation.  We do not actually discard the term involving $K$, as it will help to cancel out a higher order correction later in the derivation.  We therefore set
\begin{equation}
  \label{LM:eq:dtheta1}
  \delta \theta_1 = -(J^TJ + \lambda D^TD)^{-1} J^T r,
\end{equation}
which is the usual Levenberg-Marquardt step.

With this definition of $\delta \theta_1$, Eq.~(\ref{LM:eq:normal}) becomes
\begin{eqnarray}
  \label{LM:eq:normal2}
   \sum_{m\nu} \left( J_{m\mu} J_{m\nu} + r_m K_{m\mu\nu} + \lambda D_{m\mu} D_{m\nu} \right) \delta \theta_{2\nu} & & \nonumber \\ + \frac{1}{2} \sum_{m\nu\alpha} J_{m\mu} K_{m\nu\alpha} \delta \theta_{1\nu} \delta \theta_{1\alpha}  + \sum_{m\alpha} \left( r_m K_{m\mu\alpha}  + \delta \theta_1^\nu J_{m\nu}K_{m\mu\alpha} \right) \delta \theta_{1\alpha} & = & 0. 
\end{eqnarray}
to second order, with the term 
$r_m K_{m\mu\alpha} \delta \theta_{1\alpha}$
the term neglected by the Gauss-Newton approximation at first order.

We now turn our attention to the second term in parentheses in Eq.~(\ref{LM:eq:normal2}).  Using the definition of $\delta \theta_1 = -(J^TJ + \lambda D^TD)^{-1} J^T r$, we can write
\begin{eqnarray}
  \sum_m r_m K_{m\mu\alpha}  + \sum_{m\nu} \delta \theta_{1\nu} J_{m\nu}K_{m\mu\alpha}  & = &
  \sum_m r_m K_{m\mu\alpha}  - \sum_{m\beta\nu} r_m J_{m\beta} (J^TJ + \lambda D^TD)^{-1}_{\beta \nu} J_{n\nu} K_{n\mu\alpha} \nonumber \\
\  & = & \sum_{mn} r_m \left( \delta_{mn} - \sum_{\beta\nu} J_{m\beta} (J^TJ + \lambda D^TD)^{-1}_{\beta \nu} J_{n\nu} \right) K_{n\mu\alpha}. \nonumber
\end{eqnarray}
Since this term is proportional to the residuals, $r_m$, it can be ignored using the usual small-residual arguments.  However, we now make an appeal to geometric considerations by noting that $\delta_{mn} - \sum_{\beta\nu} J_{m\beta} (J^TJ + \lambda D^TD)^{-1}_{\beta \nu} J_{n\nu} = P^N_{mn}$ is a matrix that projects vectors perpendicular to the tangent plane of the \emph{Model Graph} as described in reference\cite{Transtrum2011}.  If the curvature of the model graph is small, then $P^N K \approx 0$ and this term can be neglected.  We discuss the implications of this argument further in section \ref{sec:smallcurvature}.

Returning to Eq.~(\ref{LM:eq:normal2}), after ignoring the last term in parentheses, we find
\begin{eqnarray}
  \label{LM:eq:dtheta2}
  \delta \theta_2 & = & -\frac{1}{2} \left( J^TJ + r^TK + \lambda D^TD \right)^{-1} J^T r'' \nonumber \\
  \ & \approx &  -1/2\, \left( J^TJ + \lambda D^TD \right)^{-1} J^T r'',
\end{eqnarray}
where we have introduced the directional second derivative $r_m'' = \sum_{\mu\nu} K_{m\mu\nu} \delta \theta_{1\mu} \delta \theta_{1\nu}$ and in the second line made the usual Gauss-Newton approximation to the Hessian, giving the formula first presented in \cite{Transtrum2010}.  This formula was originally interpreted as the second order correction to geodesic flow on the model graph, and so we refer to it as the geodesic acceleration correction.  By analogy, we refer to the first order term as the geodesic velocity.  The full step is therefore given by:
\begin{equation}
  \label{LM:eq:veldef}
  \delta \theta = \delta \theta_1 + \delta \theta_2 \equiv v \delta t + 1/2\, a \delta t^2.
\end{equation}

As has been noted previously\cite{Transtrum2010, Transtrum2011}, although the geodesic acceleration correction includes second derivative information at each step of the algorithm, its calculation is not computationally intensive.  In particular, it only requires the evaluation of a directional second derivative, which is computationally comparable to a single evaluation of the residuals.  Indeed, in the absence of an analytic expression, a finite difference estimate of the relevant second derivative can be found by a single function evaluation.  In contrast, the Jacobian evaluation at each step is comparable to $N$ function evaluations.  Particularly for large problems, the computational cost of including the geodesic acceleration is negligible compared the other elements of the algorithm.

\section{Convergence}
\label{sec:Convergence}

In order to show that geodesic acceleration does not impair the convergence guarantees of the Levenberg-Marquardt algorithm, we must make a few additional modifications.  We first note that an algorithm that selects $\Delta$ directly requires a solution of Eq.~(\ref{eq:minnolam}) given the step bound, i.e.~find the value of $\lambda$ corresponding to the step bound $\Delta$.  This so-called subproblem, can be solved accurately and efficiently for the case $K=0$, using the methods described by Mor{\'e}\cite{More1983} and Nocedel and Wright\cite{Nocedal2000}.  When including geodesic acceleration however, such a simple solution does not exist.  Indeed, the step size is no longer a monotonically decreasing function of $\lambda$.  Furthermore, accounting for the contribution from the second term requires additional function evaluations, making an accurate solution computationally expensive.

Fortunately, convergence proofs for Levenberg-Marquardt do not require that this subproblem be solved accurately.  We therefore content ourselves by approximately solving the problem as follows: We first require that $\vert \delta \theta_1 \vert \leq \Delta$.  This step can be satisfied easily using the algorithm in references\cite{More1983,Nocedal2000}.  We next require that the relative contribution from the second order step be bounded
\begin{equation}
  \label{LM:eq:va}
  \frac{ 2 \vert \delta \theta_2 \vert} {\vert \delta \theta_1 \vert} \leq \alpha
\end{equation}
for some $\alpha > 0$ \footnote{In practice $\alpha \approx 0.75$ is a robust and efficient choice}.  In practice we implement the requirement in Eq.~(\ref{LM:eq:va}) by rejecting all proposed steps for which it is not satisfied and decreasing $\Delta$ (or increasing $\lambda$) until an acceptable step is generated, recalculating $\delta \theta_1$ for the new $\Delta$.  This may appear inefficient since the method will on occasion reject steps that would have decreased the cost.  However, this requirement adds stability to the algorithm, helping it to avoid undesirable fixed points with infinite parameter values as argued and numerically justified in references\cite{Transtrum2010,Transtrum2011,Transtrum2012b}.

We need to make an additional assumption about the behavior of the model.  Specifically, it is necessary to assume that the directional second derivative of the model is bounded: $ \vert \sum_{\mu\nu} K_{m\mu\nu} u^\mu u^\nu \vert < \kappa$ for any parameter-space unit vector $u$ and some positive constant $\kappa$.  This assumption is necessary to guarantee that we can always find a step-size that satisfies Eq.~(\ref{LM:eq:va}).  We do not anticipate this requirement to be a major restriction for the applicability of the algorithm as we discuss later in this section.

With this additional assumption, we can show that our modified Levenberg-Marquardt algorithm enjoys the same global convergence properties as original algorithm.  Indeed, the proof is nearly identical to that of the unmodified Levenberg-Marquardt.  Our proof here follows closely that of Theorem 4.5 in Nocedal and Wright\cite{Nocedal2000}.  First we define the model function $m_k(\delta \theta)$ by
\begin{eqnarray}
  \label{eq:modelfunc}
  m_k(\delta \theta) & = & \frac{1}{2} \left(  r_k + J_k \delta \theta \right)^2 \\
& = & \frac{1}{2} \vert r_k \vert^2 + \delta \theta^T J_k^T r_k + \frac{1}{2} \delta \theta^T J_k^T J_k \delta \theta, \nonumber
\end{eqnarray}
and the reduction ratio $\rho_k$ by
\begin{equation}
  \label{eq:rho}
  \rho_k = \frac{ C(\theta_k) - C(\theta_k + \delta \theta_k)}{m_k(0) - m_k(\delta \theta_k)}.
\end{equation}
We now define an algorithm for which we will prove convergence.  This algorithm is analogous to algorithm 4.1 in reference\cite{Nocedal2000} which we recover if we were to set $\delta \theta_2 = 0$ at each step.

\textbf{Algorithm 1}

\begin{algorithmic}
\State Given $\hat{\Delta} > 0$, $\Delta_0 \in (0, \hat\Delta)$, and $\alpha > 0$
\For{$k = 0, 1, 2, \dots$}
   \State Calculate $\delta \theta_1$ and $\delta \theta_2$ as described above
   \State Set $\delta \theta_k = \delta \theta_1 + \delta \theta_2$
   \If{ $\vert \delta \theta_2 \vert > \alpha \vert \delta \theta_1 \vert / 2$ }
      \State $\Delta_{k+1} = \frac{1}{4} \Delta_k$
      \State $\theta_{k+1} = \theta_k$
   \Else
      \State Evaluate $\rho_k$ as in Eq.~(\ref{eq:rho})
      \If{ $\rho_k < \frac{1}{4} $ }
         \State $\Delta_{k+1} = \frac{1}{4} \Delta_k$
      \Else
         \If{ $\rho_k > \frac{3}{4}$ and $\vert \delta \theta_1 \vert = \Delta_k$ }
            \State $\Delta_{k+1} = \min(2 \Delta_k, \hat{\Delta})$
         \Else
            \State $\Delta_{k+1} = \Delta_k$
         \EndIf
      \EndIf
      \If{ $\rho_k > 0$ }
         \State $\theta_{k+1} = \theta_k + \delta \theta_k$
      \Else
         \State $\theta_{k+1} = \theta_k$
      \EndIf
   \EndIf
\EndFor 
\end{algorithmic}

Before presenting our proof, first notice that by defining $\tilde{\theta} = D \theta$, the optimization problem in $\tilde{\theta}$ must at each step satisfy the bound $\vert \delta \tilde{\theta}_1 \vert < \Delta$.  Without loss of generality, we therefore assume that $D^TD$ is the identity.  This assumption essentially replaces the Jacobian matrix with $\tilde{J} = J D^{-1}$.  With this additional assumption, we now present a Lemma that will be useful in proving convergence.

\textbf{Lemma 1}

Suppose that at some point $\theta$ our function has (non-infinite) residuals $r$, Jacobian matrix $J$, and second derivative array $K$, which satisfy $\vert J^TJ \vert < \beta$ and $ \vert K_{m\mu\nu} u^\mu u^\nu \vert < \kappa$ for any parameter space unit vector $u$.  Given positive constants $\alpha>0$ and $\zeta > 1$, then if
\begin{eqnarray}
  \label{eq:lemma1}
  \zeta \Delta & \leq & \frac{  \vert g \vert } {\sqrt{ \beta \kappa \vert r \vert / \alpha} + \beta}, \\
  \label{eq:lemma2}
  \vert \delta \theta_1 \vert & \leq & \zeta \Delta, 
\end{eqnarray}
then $\vert \delta \theta_2 \vert / \vert \delta \theta_1 \vert < \alpha/2$, where $g = J^T r$ is the function's gradient.

\textbf{Proof}

Let $\lambda$ denote the Lagrange multiplier associated with the constraint in Eq.~(\ref{eq:lemma2}).  Then $\zeta \Delta \geq \vert \delta \theta_1 \vert = \vert (J^TJ + \lambda)^{-1} g \vert \geq \vert g \vert / (\beta + \lambda)$, from which it follows that $\lambda \geq \vert g \vert / (\zeta \Delta) - \beta \geq \sqrt{\beta \kappa \vert r \vert / \alpha}$.  Notice that since $\vert J^TJ \vert \leq \beta$, the largest singular value of $J$ must be less than $\sqrt{\beta}$.  We therefore have
\begin{eqnarray}
  \vert \delta \theta_1 \vert & = & \vert \left( J^TJ + \lambda \right)^{-1} J^T r \vert \nonumber \\
  & \leq & \vert (J^TJ + \lambda)^{-1} \vert \vert J^T \vert \vert r \vert \nonumber \\
  \label{eq:dtheta1bound}
    & \leq & \frac{\sqrt{\beta} \vert r \vert}{\lambda}.
\end{eqnarray}
Similarly,
\begin{equation}
  \label{eq:dtheta2bound}
  \vert \delta \theta_2 \vert \leq \frac{\sqrt{\beta} \kappa}{2 \lambda} \vert \delta \theta_1 \vert^2.
\end{equation}
Combining these results gives us $\vert \delta \theta_2 \vert / \vert \delta \theta_1 \vert \leq \beta \kappa \vert r \vert / 2 \lambda^2< \alpha/2$.

With this Lemma, we are prepared to prove our main result: that including geodesic acceleration does not affect the convergence properties of the Levenberg-Marquardt algorithm.

\textbf{Theorem 1}

Suppose that $\vert J^T J \vert \leq \beta$ for some positive constant $\beta$, that $C$ is Lipschitz continuously differentiable in the neighborhood $S(R_0)$ for some $R_0 > 0$, and that $\vert \delta \theta_1 \vert < \zeta \Delta_k$ for some constant $\zeta > 1$ at each iteration.  Also assume that $m_k(0) - m_k(\delta \theta) \geq c_1 \vert g_k \vert \min( \Delta_k, \vert g_k \vert / \vert J^T J \vert )$ for some positive constant $c_1 \in (0, 1]$ ($g_k = J^T_k r_k$), and that the second directional derivative of $r_k$ in any unit direction $u$ is bounded $\vert \sum_{\mu\nu} K_{m\mu\nu} u^\mu u^\nu \vert < \kappa$ for some positive $\kappa$.  Then using Algorithm 1,
\begin{equation}
  \label{eq:themresult}
  \lim_{k\rightarrow \infty} \inf \vert g_k \vert = 0.
\end{equation}

\textbf{Proof}

The proof is nearly identical to that of Theorem 4.5 in reference \cite{Nocedal2000}.  We repeat the proof here in order to highlight the small-differences, leaving out the algebraic details.

With some algebraic manipulation, we obtain
\[ \vert \rho_k - 1 \vert = \vert \frac{ m_k(\delta \theta_k) - C(\theta_k  + \delta \theta_k)}{ m_k(0) - m_k(\delta \theta_k)} \vert. \]
Using the same argument as in reference \cite{Nocedal2000}, we have
\begin{equation}
  \label{eq:427}
  \vert m_k(\delta \theta_k) - C(\theta_k + \delta \theta_k) \vert \leq (\beta/2 + \beta_1) \vert \delta \theta_k \vert^2, 
\end{equation}
where $\beta_1$ is the Lipschitz constant for $g$ on the set $S(R_0)$.

Suppose for contradiction that there is $\epsilon > 0$ and a positive index $K$ such that
\begin{equation}
  \label{eq:428}
  \vert g_k \vert \geq \epsilon, \ \ \ \textrm{for all} \ k \geq K,
\end{equation}
then we have
\begin{equation}
  \label{eq:429}
  m_k(0) - m_k(\delta \theta_k) \geq c_1 \vert g_k \vert \min(\Delta_k, \vert g_k \vert / \vert J^T J \vert ) 
  \geq c_1 \epsilon \min(\Delta_k, \epsilon/\beta).
\end{equation}

By the workings of Algorithm 1, we see that $\rho$ is only calculated if $\vert \delta \theta_2 \vert < \alpha \vert \delta \theta_2 \vert /2 $.  In this case we have $\vert \delta \theta \vert = \vert \delta \theta_1 + \delta \theta_2 \vert < \vert \delta \theta_1 \vert + \vert \delta \theta_2 \vert < (1 + \alpha/2) \vert \delta \theta_1 \vert < (1 + \alpha/2) \zeta \Delta_k$.  Each step therefore satisfies
\begin{equation}
  \label{eq:425}
  \vert \delta \theta \vert < \gamma \Delta_k
\end{equation}
where $\gamma = (1 + \alpha/2) \zeta > 1$.

Using Eqs.~(\ref{eq:425}), (\ref{eq:429}), and (\ref{eq:427}), we have
\begin{equation}
  \label{eq:430}
  \vert \rho_k - 1 \vert \leq \frac{ \gamma^2 \Delta_k^2 (\beta/2 + \beta_1)}{c_1 \epsilon \min(\Delta_k, \epsilon/\beta)}.
\end{equation}
We now derive a bound on the right-hand-side that holds for all sufficiently small-values of $\Delta_k$, that is, for all $\Delta_k \leq \bar{\Delta}$, where $\bar{\Delta}$ is defined as follows:
\begin{equation}
  \label{eq:431}
  \bar{\Delta} = \min \left(  \frac{1}{2} \frac{c_1 \epsilon}{\gamma^2 (\beta/2 + \beta_1)} ,\frac{R_0}{\gamma}, \frac{\epsilon}{ \zeta \sqrt{ \beta \kappa \vert r \vert / \alpha} + \zeta \beta} \right).
\end{equation}
As noted in reference \cite{Nocedal2000}, the $R_0/\gamma$ term in this definition ensures that the bound in Eq.~(\ref{eq:427}) is valid because $\vert \delta \theta_k \vert \leq \gamma \Delta_k \leq \gamma \bar{\Delta} \leq R_0$.  The third term is necessary to show convergence when including geodesic acceleration.  

Note that since $c_1 \leq 1$ and $\gamma \geq 1$, we have $\bar{\Delta} \leq \epsilon/\beta$.  The latter condition implies that for all $\Delta_k \in [0, \bar{\Delta}]$, we have $\min(\Delta_k, \epsilon/\beta) = \Delta_k$, so from Eq.~(\ref{eq:430}) and (\ref{eq:431}), we have 
\[  \vert \rho_k - 1 \vert \leq \frac{1}{2}, \]
following the logic in reference \cite{Nocedal2000}.  Therefore, if $\Delta_k \in [0, \bar{\Delta}]$, then $\rho_k > \frac{1}{4}$.  Furthermore, the third condition in Eq.~(\ref{eq:431}) implies that if $\Delta_k \in [0, \bar{\Delta}]$ then $\vert \delta \theta_2 \vert < \alpha \vert \delta \theta_1 \vert /2$ by Lemma 1.

Since, we have that if $\Delta_k \in [0, \bar{\Delta}]$, then $\vert \delta \theta_2 \vert < \alpha \vert \delta \theta_1 \vert /2$ and that $\rho_k > \frac{1}{4}$, by the workings of Algorithm 1, $\Delta_{k+1} \geq \Delta_k$ whenever $\Delta_k$ falls below the threshold $\bar{\Delta}$.  Consequently, a reduction in $\Delta_k$ (by a factor of $1/4$) can occur in the algorithm only if $\Delta_k \geq \bar{\Delta}$.  We conclude that
\begin{equation}
  \label{eq:432}
  \Delta_{k+1} \geq \min( \Delta_k, \bar{\Delta}/4) \ \ \ \textrm{for all} \ k \geq K.
\end{equation}

Suppose now that there is an infinite subsequence $\mathcal{K}$ such that $\rho_k \geq 1/4$ for $k\in\mathcal{K}$.  For $k \in \mathcal{K}$ and $k \geq K$, we have from Eq.~(\ref{eq:429}) that
\begin{eqnarray}
  C(\theta_k) - C(\theta_{k+1}) & = & f(\theta_k) - f(\theta_k + \delta \theta_k) \nonumber \\
  & \geq & \frac{1}{4} [ m_k(0) - m_k(\delta \theta_k) ] \nonumber \\
  & \geq & \frac{1}{4} c_1 \epsilon \min( \Delta_k, \epsilon/\beta). \nonumber
\end{eqnarray}
Since $C$ is bounded below, it follows from this inequality that 
\[ \lim_{k\in\mathcal{K}, k \rightarrow \infty} \Delta_k = 0, \]
contradicting Eq.~(\ref{eq:432}).  Hence no such infinite subsequence $\mathcal{K}$ can exist, and we must have $\rho_k < \frac{1}{4}$ for all $k$ sufficiently large.  In this case, $\Delta_k$ will eventually be multiplied by $\frac{1}{4}$ at every iteration, and we have $\lim_{k\rightarrow \infty} \Delta_k = 0$, which again contradicts Eq.~(\ref{eq:432}).  Hence our original assertion, Eq.~(\ref{eq:428}) must be false, proving the theorem.

Although Theorem 1 and its subsequent proof are nearly identical to Theorem 4.5 in Nocedal and Wright\cite{Nocedal2000}, we now briefly discuss its implication.  In particular, the proof is that the algorithm converges to a point where cost has zero gradient.  This does not automatically imply convergence in the parameters unless the Hessian matrix is bounded from below.  In fact, the Hessian matrix can be very ill-conditioned, particularly for models with many parameters\cite{Waterfall2006,Gutenkunst2007a}.  Although the inferred parameters in such problems may be ill-conditioned, the convergence properties of the algorithm are nevertheless robust.

The requirement that $m_k(0) - m_k(\delta \theta) \geq c_1 \vert g_k \vert \min( \Delta_k, \vert g_k \vert / \vert J^T J \vert )$ at each step of the algorithm in this context is given without motivation.  However, just like the analogous theorem in reference\cite{Nocedal2000}, this requirement is necessary to guarantee that the parameter iterates do not accumulate at points with nonzero gradient.  We note that in practice this condition is never explicitly checked.  

The only new assumption of Theorem 1 beyond those for the standard algorithms\cite{Nocedal2000} is the bound on the second derivative. Since we are utilizing second derivative information, this addition is not unexpected and guarantees that the new algorithm is well-behaved. It is also clearly analogous to the bound on the first derivative used in the first algorithm ($|J^TJ| <= \beta$, equivalent to $|J u| <= \sqrt{\beta}$). The standard bound on the first derivative already excludes models where the cost diverges at interior points in the domain of the model as well as points like $C(\theta) \sim \sqrt{|\theta-\theta_0|}$ where the cost would stay finite although the derivative diverges; the new bound on the second derivative additionally excludes cases like $C(\theta) \sim |\theta - \theta_0|^{3/2}$, which should not arise often in practice. In our experience, many models do have costs which become singular at unphysical points ({\it i.e.}, positive-only parameters set less than zero), which can be addressed by an appropriate reparameterization of the model ({\it i.e.}, shifting to log parameters). When this can be done, it is also likely to improve the convergence rate of the algorithm.

We have here shown convergence of geodesic acceleration for an algorithm that belongs to the broad class of trust region methods that operate by specifying a step bound $\Delta$.  Originally Levenberg-Marquardt was proposed as an algorithm that heuristically selected $\lambda$ directly rather than implicitly through $\Delta$\cite{Levenberg1944,Marquardt1963}.  For an algorithm that directly selects $\lambda$, one can similarly show convergence\cite{Osborne1976}.  The proof that geodesic acceleration does not impair the convergence of this class of algorithm is almost identical to that of the original theorem, just as the proof of Theorem 1 above is almost identical to that of Theorem 4.5 in reference\cite{Nocedal2000}.  We do not give the rigorous proof here, but merely note that under the same additional assumptions, convergence can be shown for these methods as well.

\section{Small-curvature approximation and relation to other methods}
\label{sec:smallcurvature}

In deriving the geodesic acceleration in section \ref{sec:derivation}, we observed that the neglected terms were each proportional to the residuals and could be neglected in the small-residual approximation.  We also noted that the neglected terms were also proportional to the extrinsic curvature on the model graph as described in reference\cite{Transtrum2010,Transtrum2011}.  Indeed, the geodesic acceleration was originally understood as a small-curvature approximation rather than a small-residual approximation.  These two approximations are complementary; in general only one of the two needs to  be valid in order for the approximations to hold.  

In deriving the geodesic acceleration, the connection between the small-residual and small-curvature approximation only became apparent when considering cubic order terms.  However, the equivalence of the two approximations can be seen in the neglected terms of the Gauss-Newton Hessian without considering the geodesic acceleration.  At a fixed point of the cost, the residuals are perpendicular to the model manifold (a fact used elsewhere to justify a scale-free measure of convergence\cite{Bates1981b}), i.e. $r_n \approx \sum_m r_{m} P^N_{mn}$.  Thus, as the algorithm approaches a minimum, the neglected term can be written as
\begin{equation}
  \label{eq:curvaturehessian}
  \sum_m r_m \frac{\partial^2 r_m}{\partial \theta_\mu \partial \theta_\nu}  \approx \sum_{mn} r_m P^N_{mn} \frac{\partial^2 r_n}{\partial \theta_\mu \partial \theta_\nu}.
\end{equation}
As Eq.~(\ref{eq:curvaturehessian}) makes clear, near the best fit the nonlinear contribution to the Hessian includes only components perpendicular to the tangent plane and is thus proportional to the extrinsic curvature.

The implication of Eq.~(\ref{eq:curvaturehessian}) is that the Levenberg-Marquardt algorithm may attain very good convergence rates even with large residuals so long as the extrinsic curvature is sufficiently small.  However, the approximation in Eq.~(\ref{eq:curvaturehessian}) is only valid near a fixed point of the cost, just as is the small-residual approximation is only valid near the minimum.  The advantage of identifying the equivalence of the small-residual and small-curvature approximations is that the latter is likely to have much wider applicability.  In particular for data fitting, the small-curvature approximation is a feature of the model rather than the data.  Thus, although the validity of the small-residual approximation cannot be identified without finding the best fit, the small-curvature approximation can be.  Furthermore, the small-curvature approximation is likely to be an excellent approximation for most models.  In particular, several examples in references\cite{Bates1980,Transtrum2011} exhibit extrinsic curvatures many orders of magnitude smaller than the magnitude of the bare nonlinearities.  

Although Eq.~(\ref{eq:curvaturehessian}) is only valid near a minimum, the approximation that led to the geodesic acceleration is valid over a much larger parameter range, and it is for this reason it is likely to be a useful modification.  The problem with Levenberg-Marquardt is not that its asymptotic convergence rate is poor (super-linear convergence is satisfactory in most cases), but that it may spend an unreasonable amount of time navigating a narrow canyon before that convergence rate is realized.  As originally suggested, the geodesic acceleration can speed up this process since it describes the curvature of the canyon.  The algorithm can find the minimum more quickly by following the canyon with a sequence of parabolic steps, rather than linear steps.  A numerical demonstration of this will not be given here, as it has already been shown elsewhere\cite{Transtrum2010,Transtrum2011,Transtrum2012b}.

There are many other methods which include second derivative information in order to improve the Levenberg-Marquardt algorithm.  These approaches all try to improve the quality of the Hessian estimate in some way\cite{DennisJr1981,Fletcher1987,Nocedal2000}.  Thus, although geodesic acceleration may appear superficially similar to these approaches, it is actually quite different.  The philosophy behind geodesic acceleration is to extend the small-residual/curvature approximation to higher order terms rather than estimate the neglected terms.  Of course, the utility of such of an approach is ultimately measured by performance on real world problems.  As has been shown elsewhere, geodesic acceleration can be very helpful on large problems for which Levenberg-Marquardt spends unreasonable time searching for a minimum before zooming into the best fit.

\section{Conclusion}
\label{sec:conclusion}

In this paper we have studied the geodesic acceleration correction to the Levenberg-Marquardt algorithm, originally suggested in references\cite{Transtrum2010,Transtrum2011}.  We have derived the correction without the use of differential geometry and shown that it can be interpreted as an extension of the small-residual approximation used to estimate the Hessian matrix of a sum of squares.  We have also seen that the small-residual approximation is complemented by the small-curvature approximation, which is likely to be applicable under much more general circumstances.  

We have seen that with just a few minor modifications, the geodesic acceleration correction can be incorporated into a standard Levenberg-Marquardt routine with minimal computational cost and without sacrificing its convergence properties.  Numerical experiments given elsewhere\cite{Transtrum2010,Transtrum2011,Transtrum2012b}, suggest that the benefit/cost ratio of this improvement can be substantial on many optimization problems.

The authors would like to thank Cyrus Umrigar, Peter Nightingale, Saul Teukolsky, and Ben Machta for helpful conversation.  This work is partially supported by the TCGA Genome Data Analysis Center (GDAC) grant (MKT) and the Cancer Center Support Grant at the University of Texas MD Anderson Cancer Center (U24 CA143883 02 S1 and P30 CA016672) (MKT) and by NSF grant number DMR-1005479 (JPS).

\bibliographystyle{aps}
\bibliography{C:/Users/MKTranstrum/Z/References}

\end{document}